\documentstyle[12pt,amstex,amssymb]{article}
\renewcommand{\a}{\alpha}

\renewcommand{\d}{\delta}
\newcommand{\D}{\Delta}

\renewcommand{\l}{\lambda}
\renewcommand{\L}{\Lambda}

\newcommand{\s}{\sigma}

\renewcommand{\i}{\infty}
\newcommand{\p}{\partial}

\title{The Thermal Explosion Revisited}
\author{G. I. Barenblatt,${}^{\dag,\S}$ J. B. Bell,${}^{\S}$ W. Y.
Crutchfield${}^{\S}$}
\date{}

\begin{document}

\maketitle

\begin{center}
Contributed by Grigory Isaakovich Barenblatt
\end{center}

\vfill

\noindent
\underline{\hskip 1 in}

\noindent
${}^{\dag}$Department of Mathematics, University of California at Berkeley,
Berkeley, CA\ \ 94720

\noindent
${}^{\S}$Lawrence Berkeley National Laboratory, Berkeley, CA\ \ 94720

\newpage
\begin{abstract}
{\bf The classical problem of the thermal explosion in a long cylindrical
vessel is modified so that only a fraction $\a$ of its wall is ideally
thermally conducting while the remaining fraction $1-\a$ is thermally
isolated.  Partial isolation of the wall naturally reduces the critical
radius of the vessel.  Most interesting is the case when the structure of
the boundary is a periodic one, so that the alternating conductive $\a$ and
isolated $1-\a$ parts of the boundary occupy together the segments $2\pi/N$
($N$ is the number of segments) of the boundary.  A numerical investigation
is performed.  It is shown that at small $\a$ and large $N$ the critical
radius obeys a scaling law with the coefficients depending upon $N$.  For
large $N$ is obtained that in the central core of the vessel the temperature
distribution is axisymmetric.  In the boundary layer near the wall having
the thickness $\approx 2\pi r_0/N$ ($r_0$--the radius of the vessel) the
temperature distribution varies sharply in the peripheral direction.  The
temperature distribution in the axisymmetric core at the critical value of
the vessel radius is subcritical.}
\end{abstract}

\newpage
\section{Introduction}

We revisit in this Note the classical problem of a thermal explosion in a
long circular cylindrical vessel containing an exothermically reacting gas
at rest.  It was formulated and solved under some assumptions by D.~A.
Frank--Kamenetsky \cite{1}, \cite{2} (see also \cite{3}).  In the original
formulation it was assumed that the wall of the vessel is ideally
conducting, so that the gas temperature at the boundary is equal to the
temperature of the ambient medium.  Such boundary condition made the
problem cylindrically symmetric, and the symmetry essentially simplified
the solution.  In the present Note the problem is modified in the
following way:
the wall
is partially isolated so that the symmetry is lost.  The critical values of
the radius of the vessel are determined numerically.  Especially
instructive results were obtained for the cases when the isolated parts of
the boundary are distributed periodically with large angular frequency, and
the isolated part of the boundary is large.  The scaling laws for the
critical values were found.  For large angular frequencies it was found that
there exists an axisymmetric core of the temperature distribution which
occupies a major part of the vessel.  The conditions in this core at the
critical case were found to be subcritical.

\section{Mathematical Problem Formulation}

Assume that a gas at rest is enclosed in a long cylindrical vessel of radius
$r_0$.  An exothermic reaction is going in the gas with the thermal effect
$Q$ per unit mass of reacted gas.  For the reaction rate the Arrhenius law is
assumed with the activation energy $E$.  If the thermal effect and the
activation energy are large, it can be shown (see \cite{3}) that an
`intermediate-asymptotic' steady state regime is achieved.  For this regime
the gas consumption in the reaction can be neglected, and, from the other
side, the temperature distribution in the vessel is steady.  Applying the
Frank--Kamenetsky large activation energy approximation, a non-linear
equation for dimensionless reduced temperature ${\bold u}$ is obtained:
\setcounter{equation}{0}
\begin{equation}
\label{eq2.1}
\D u + \l^2 e^u = 0,
\end{equation}
where
\begin{equation}
\label{eq2.2}
u = \frac {(T-T_0)E}{RT_0^2},
\end{equation}
the Laplace operator $\D$ is related to dimensionless variables $\rho =
r/r_0$, $\theta$; $r,\theta$ are the polar coordinates.  The constant $\l$ is
\begin{equation}
\label{eq2.3}
\l = \frac {r_0}{l},\ l = \left( e^{\frac
{E}{RT_0}}\kappa RT_0^2c/QE\sigma(T_0)\right)^{1/2}.
\end{equation}
Here $T$ is the absolute temperature, $T_0$--the temperature of the ambient
medium, $R$--the universal gas constant, $\kappa$--the molecular thermal
diffusivity, ${\bold c}$--heat capacity of gas per unit volume, $\s(T)$ is
the pre-exponential factor in the Arrhenius reaction rate expression:  a
slow function of temperature.

In the classical problem formulation it was assumed that the whole wall of
the vessel is ideally heat conducting, so that the gas temperature at the
boundary is equal to the temperature of the ambient medium.  This gives a
Dirichlet condition for the equation (\ref{eq2.1}):
\begin{equation}
\label{eq2.4}
u(1,\theta) = 0.
\end{equation}

The boundary value problem formulation under this condition is axisymmetric,
and this was important for obtaining the analytic solution in an explicit
form.  D.~A. Frank--Kamenetsky showed that the solution to the problem
(\ref{eq2.1}), (\ref{eq2.4}) does exist for $\l \le \l_{cr} = \sqrt{2}$
only.  Physically it means that a quiet, non-explosive proceeding of the
reaction is possible only if the radius of the vessel is less than a
critical one:  $r_0 \le (r_0)_{cr} = \sqrt{2} l$.  This condition is known
as the {\em condition of the thermal explosion}.

In the present Note the following modification of the problem (\ref{eq2.1}),
(\ref{eq2.4}) is proposed:  Only a fraction $\a$ of the wall is heat
conducting, while the fraction $1-\a$ is thermally isolated.  The simplest
formulation corresponds to a mixed problem (Figure 1,a)
\begin{equation}
\label{eq2.5}
\begin{array}{rll}
u(1,\theta) = 0 &\mbox{at} &0 \le \theta \le 2\pi \a \vspace{1\jot} \\
\p_{\rho}u(1,\theta) = 0 &\mbox{at} &2\pi \a < \theta \le 2\pi.
\end{array}
\end{equation}

More interesting is the case when the isolated part of the wall is not
concentrated on a single arc, but is distributed periodically (Figure 1,b):
the boundary $\rho = 1$ is divided into $N$ segments
\begin{equation}
\label{eq2.6}
\frac {2\pi}{N} (1+k) \ge \theta \ge \frac {2\pi}{N} k,\quad k =
0,1,\dots,N-1.
\end{equation}

The fraction $\a$ of each segment is left heat conducting, while the fraction
$1-\a$ becomes isolated.  In this case the mixed boundary condition at $\rho
= 1$ has the form:
\begin{equation}
\label{eq2.7}
\begin{array}{rll}
u(1,\theta) = 0 &\mbox{at} &\displaystyle{\frac {2\pi}{N} k \le \theta \le
\frac {2\pi \a}{N} + \frac {2\pi}{N} k}, \vspace{1\jot} \\
\p_{\rho}u(1,\theta) = 0 &\mbox{at} &\displaystyle{\frac {2\pi\a}{N} + \frac
{2\pi}{N} k < \theta \le \frac {2\pi}{N} (k+1)}.
\end{array}
\end{equation}

The central question addressed in the present Note is:  What are the
asymptotic laws for the critical radius if $N \to \i$ and $\a \to 0$, i.e.
the period of the boundary condition tends to zero, and the isolation is
close to a complete one.  Remember:  for the complete isolation the critical
radius is equal to zero.

\section{The Numerical Method}

In order to answer the questions posed above, we must numerically evaluate
the critical values $\l_{cr}$ for fixed $N$ and $\a$.  There are two
special items in determining $\l_{cr}$: the singularity of the
linearized equation (\ref{eq2.1}):
\setcounter{equation}{0}
\begin{equation}
\label{eq3.1}
\D \d u + \l^2e^u \d u = 0
\end{equation}
at the critical value of $\l = \l_{cr}$, and the dependence of the values of
$\l_{cr}$ obtained by discretization upon the number of
points, $n$ used to discretize each dimension.  Here $\d u$ is the
perturbation of the solution.

For fixed values of $N$, $n$, and $\a$, we determined a trajectory of
solutions $u$ versus $\l$ by solving equation (\ref{eq2.1}) with a
Newton--Raphson method.  This method requires the solution of the linearized
equation (\ref{eq3.1}) which becomes singular at $\l = \l_{cr}$.  In
practice it prevents us from approaching the critical point closely.
Therefore to make an accurate determination of $\l_{cr}$ an extrapolation
procedure was used.  It was assumed that for $\l$
approaching $\l_{cr}$ a parabolic approximation is valid:
\begin{equation}
\label{eq3.2}
\l_{cr}^2 - \l^2 = C(\|u\| - u_0)^2.
\end{equation}
The parameters $\l_{cr}$, $C$ and $u_0$ were determined to fit the last
10 points on the trajectory $u$ versus $\l$ approaching $\l_{cr}$.  The
approximation (\ref{eq3.2}) happened to be satisfactory.  Typically the
fit (\ref{eq3.2}) is accurate to a few parts in $10^6$.

The above procedure
yielded an estimate for
$\l_{cr}$ as a function of
$N$, $\a$ and the number of discretization points $n$.  It is natural to
remove the dependence of $\l_{cr}$ of non-physical, computational parameter
$n$.  For this purpose another extrapolation was used.  In our numerical
approximations the second-order accurate discretizations of the operators
was employed.  If the boundary conditions were smooth, the solution $u$ would
approach a limit with an error of the order of $1/n^2$.
But the boundary conditions are non-smooth, the derivative of the solution is
discontinuous at the boundary, and it causes the order of the approximation
to decrease.  Extensive
numerical calculations have shown $\l_{cr}^2$ to vary linearly with
$1/n$.  Therefore, the following iterative procedure was used:  the value
$\l_{cr}(N,\a,n)$ in (\ref{eq3.2}) was calculated for three different values
of $n$.  These three values are then fit to a linear function $a+b/n$.  If
the fit was poor, as might happen if the values of $n$ were too
small, the procedure is repeated with larger values of $n$ and so
on, until a satisfactory fitting was obtained.  The extrapolation $n \to \i$
is simply the value of the coefficient $a$.

\section{Results of the Numerical Analysis}

The results obtained by numerical solution are represented on Figures 2--4.
Three instructive properties are revealed.

(i) On the graph of Figure 2 the values of $\l_{cr}^2$ are presented for
growing values of $N$ as the functions of $1/\a$.  It is seen that the
critical value $\l_{cr}$, i.e. the critical radius of the vessel for large
$N$ is practically insensitive to $\a$ up to $\a$ very small.  For small $N$
the dependence of $\l_{cr}$ on $\a$ is strong.  Clearly, for any $N$,
$\l_{cr} = 0$ for $\a = 0$, but it is instructive that for instance, for $N
= 256$ when only $1/512$ ($0.2$ percent) part of the boundary is heat
conducting, the
critical value of the radius is only 4 percent less than the critical radius
for wholly heat conducting wall.

(ii) For large $N$, starting, say, from $N = 32$, there exists an internal
core $0 \le \rho \le \rho_*$ where the solution is close to axisymmetric one
(see Figure 3).  The value $\rho_*$ was selected so that
$|u_{\max}(\rho_*,\theta) - u_{\min}(\rho_*,\theta)| < 10^{-4}$.
Introducing the mean value $u_* = (u_{\max} + u_{\min})/2$ we notice, that
for $0 \le \rho \le \rho_*$ the solution is close to axisymmetric, so that
the equation (\ref{eq2.1}) and the boundary condition at $\rho = \rho_*$
assume the form:
\setcounter{equation}{0}
\begin{equation}
\label{eq4.1}
\frac {1}{\rho} \ \frac {d}{d\rho} \rho \  \frac {du}{d\rho} + \l^2e^u = 0,\
u = u_* \mbox{ at } \rho = \rho_*.
\end{equation}

Transforming the variables
\begin{equation}
\label{eq4.2}
u = u_* + v,\  \ \ R = \frac {\rho}{\rho_*}
\end{equation}
we reduce the problem to a classic one:
\begin{equation}
\label{eq4.3}
\frac {1}{R} \  \frac {d}{dR} R  \ \frac {dv}{dR} + \L^2e^v = 0,\ \ \
v(1,\theta) = 0
\end{equation}
where $0 \le R \le 1$, $\L^2 = \l^2\rho_*^2e^{u_*}$.  We calculate now the
values of $\L_{cr}^2 = \l_{cr}^2\rho_*^2e^{u_*}$ where $\l_{cr}$ is the
critical value obtained in previous calculations.  The graphs $\L_{cr}^2$ as
the functions of $1/\a$ for different $N$ are presented on Figure 4.  It can
be seen that the values of $\L_{cr}^2$ are always less than $2$ (within the
limits of our numerical accuracy).  This means that the `rugged' boundary
layer near the wall $\rho = 1$ controls the approach to criticality.  The
thickness of this layer is of the order of the length of the segment
$2\pi/N$.  The angular derivative
$\p_{\theta}u$ in the boundary layer is large.  The value $u_*$ decreases as
$N$ increases.

(iii) The intermediate power laws are observed for $\l_{cr}^2$ at large $N$
and small $\a$:
\begin{equation}
\label{eq4.4}
\l_{cr}^2 = S(N)\a^{t(N)}.
\end{equation}

The values of $S(N)$ and $t(N)$ for various $N$ are given in the Table.

\section{Conclusion}

Non-axisymmetric modification of the problem of  thermal explosion in a
cylindrical vessel is formulated.  The boundary is partly isolated, and only
partly ideally conducting.  Special attention is paid to the case of
periodic distribution of the isolated and conducting parts.  The critical
values of radius and other relevant properties are obtained numerically.

It is shown that for the period small in comparison with the vessel radius
the critical value of radius of the vessel is practically insensitive to the
relative size of the open area of the wall up to its very small values.  The
temperature distribution in the central core is axisymmetric and subcritical
even at globally critical conditions:  the criticality is due to a thin
boundary layer near the wall where the temperature distribution is highly
non-axisymmetric.  Intermediate power laws are obtained for the critical
radius as the function of the relative open area of the wall.

This work was supported in part by the Applied Mathematical Sciences
subprogram of the Office of Energy Research, U.S. Department of Energy under
Contract DE-AC-03-76-SF00098, and in part by the National Science Foundation
under Grants DMS94-14631 and DMS-2732710.

\clearpage

\begin{center}
{\bf Table}
\end{center}

\begin{center}
\begin{tabular}{rll}
$N$\  &$S(N)$ &$t(N)$ \\
32 &2.03 &0.10 \\
64 &2.03 &0.055 \\
128 &2.02 &0.028 \\
256 &2.01 &0.015
\end{tabular}
\end{center}

\clearpage

\begin{center}
{\bf Figure Captions}
\end{center}

\medskip
\noindent
Figure 1.  A fraction of the wall is isolated.  (a) Isolated (5/6) and
conducting (1/6) parts are connected.  (b) Isolated and conducting parts are
distributed periodically.

\medskip
\noindent
Figure 2.  Dimensionless critical radius as the function of conducting
fraction $\alpha$ for different angular frequencies.  It is seen that at
large frequencies the critical radius is practically $\alpha$-independent up
to very small values of the conducting fraction $\alpha$.

\medskip
\noindent
Figure 3.  The solution reveals an axisymmetric internal central core and
`rugged' boundary layer ($N=32$, $\alpha = 1/32$).

\medskip
\noindent
Figure 4.  The temperature distribution in the internal central core is a
subcritical one:  $\Lambda^2 < 2$.

\end{document}